\documentclass{article}
\usepackage[utf8]{inputenc}
\usepackage[english]{babel}
\usepackage{amsmath}
\usepackage{amsthm}
\usepackage{amsfonts}
\usepackage{amssymb}
\usepackage{yfonts}

\begin{document}

\author{Walter Wyss}
\title{Two Non-Commutative Binomial Theorems}
\date{}
\maketitle

\begin{abstract}
We derive two formulae for $(A+B)^n$, where $A$ and $B$ are elements in a non-commutative, associative algebra  with identity.
\end{abstract}

\section{Introduction}
Let \textswab{A} be an associative algebra, not necessarily commutative, with identity. For two elements $A$ and $B$ in \textswab{A}, that commute, i.e. \\

\begin{equation}
AB=BA
\end{equation}

the well-known Binomial Theorem reads 
\begin{equation}
(A+B)^n = \sum_{k=0}^n \begin{pmatrix} n \\ k \end{pmatrix} A^kB^{n-k}
\end{equation}

If $A$ and $B$ do not commute, we find the first formula for $(A+B)^n$ that retains the binomial coefficient. It also gives a representation of $e^{(A+B)}$ that is different from the Campell-Baker-Hausdorff representation \cite{3}. The first formula is then applied to a problem in non-commutative geometry. The second formula for $(A+B)^n$ complements the first one. We apply it to a problem in quantum mechanics.
 
\section{The First Non-Commutative Binomial Theorem}

Let \textswab{A} be an associative algebra, not necessarily commutative, with identity 1. $L$(\textswab{A}) denotes the algebra of linear transformations from \textswab{A} to \textswab{A}.

\subsection*{Definition 1}
Let $A$ and $X$ be elements of \textswab{A}.
 \begin{enumerate}
 \item $A$ \text{can be looked upon as an element in} $L$(\textswab{A}) \text{by} \\
 
 \begin{equation}
 A(X)= AX
 \end{equation}
 
 \text{i.e. leftmultiplication}

\item The element $d_A$ in $L$(\textswab{A}) is defined by
 \begin{equation}
 d_A(X) = [A,X] = AX - XA
 \end{equation}

\text{We now have the following trivial relations:}

 \end{enumerate}

\subsection*{Statements}
\begin{enumerate}
\item As elements in $L$(\textswab{A}), $A$ and $d_A$ commute, i.e.
\begin{equation}
Ad_A(X) = d_AA(X)
\end{equation}

\item $d_A$ is a derivation on \textswab{A}, i.e.
\begin{equation}
d_A(XY)=(d_AX)Y+X(d_AY)
\end{equation}

\item
\begin{equation}
(A-d_A)X = XA
\end{equation}

\item Jacobi identity
\begin{equation}
d_Ad_B(C) + d_Bd_C(A) +d_Cd_A(B) = 0
\end{equation}

\end{enumerate}

These simple statements are sufficient to prove the following non-commutative Binomial Theorem \cite{1}, \cite{2}.

\section*{Theorem 1}
For A and B elements in \textswab{A}, and 1 being the identity in \textswab{A}

\begin{equation}
(A+B)^n = \sum_{k=0}^n \begin{pmatrix} n \\ k \end{pmatrix} \lbrace (A+d_B)^k1\rbrace B^{n-k}
\end{equation}

\begin{proof}
The formula holds true for n=1. We now proceed by induction.

\begin{align*}
(A+B)^{n+1} = (A+B) (A+B)^n = (A+d_B+B-d_B)(A+B)^n \\
= (A+d_B+B-d_B) \sum_{k=0}^n \begin{pmatrix} n \\ k \end{pmatrix} \lbrace (A+d_B)^k1\rbrace B^{n-k}
\end{align*}

Using the previous Statements, we get 

\begin{align*}
(A+B)^{n+1} = &\sum_{k=0}^n \begin{pmatrix} n \\ k \end{pmatrix} [A \lbrace (A+d_B)^k1 \rbrace B^{n-k} + \lbrace d_B(A +d_B)^k1 \rbrace B^{n-k} + \lbrace (A +d_B)^k1 \rbrace B^{n-k+1}] \\
= &\sum_{k=0}^n \begin{pmatrix} n \\ k \end{pmatrix} [\lbrace (A+d_B)^{k+1}1 \rbrace B^{n-k} + \lbrace (A+d_B)^k1 \rbrace B^{n-k+1} ] \\
= &\sum_{k=1}^n \begin{pmatrix} n \\ k \end{pmatrix} \lbrace (A+d_B)^{k}1 \rbrace B^{n-k+1} + B^{n+1} \\
& + \sum_{k=1}^n \begin{pmatrix} n \\ k-1 \end{pmatrix} \lbrace (A+d_B)^{k}1 \rbrace B^{n-k+1} + \lbrace (A+d_B)^{n+1}1 \rbrace 
\end{align*}

From the identity 
\begin{equation*}
\begin{pmatrix} n \\ k \end{pmatrix} + \begin{pmatrix} n \\ k-1 \end{pmatrix} = \begin{pmatrix} n+1 \\ k \end{pmatrix}
\end{equation*}

we then get
\begin{equation*}
(A+B)^{n+1} = \sum_{k=0}^{n+1} \begin{pmatrix} n+1 \\ k \end{pmatrix} \lbrace (A+d_B)^{k}1 \rbrace B^{n+1-k}
\end{equation*}

\end{proof}

\section{The Essential Non-Commutative Part}

We write 
\begin{equation}
(A+d_B)^n1 = A^n+D_n(B,A)
\end{equation}

For a commutative algebra, $D_n(B,A)$ is identically zero. We thus call $D_n(B,A)$ the essential non-commutative part. \\
$D_n(B,A)$ satisfies the following recurrence relation

\begin{equation}
D_{n+1}(B,A) = d_BA^n + (A+d_B)D_n(B,A)
\end{equation}

with 
\begin{equation*}
D_0(B,A) = 0
\end{equation*}

\subsection*{Definition 2}
\begin{enumerate}
\item
\begin{equation}
M_n = \sum_{k=0}^n \begin{pmatrix} n \\ k \end{pmatrix} A^kB^{n-k}
\end{equation}
\item
\begin{equation}
D_k(B,A) =D_k
\end{equation}
\end{enumerate}

We now have the following obvious corollary.

\subsection*{Corollary 1}
\begin{equation}
(A+B)^n = M_n + \sum_{k=0}^n \begin{pmatrix} n \\ k \end{pmatrix} D_k B^{n-k}
\end{equation}

\section{Exponentials}
We have as a consequence of the first non-commutative Binomial Theorem

\subsection*{Corollary 2}
\begin{equation}
e^{A+B} = [e^{A+d_B} 1]e^B 
\end{equation}
\begin{proof}
\begin{align}
e^{A+B} &= \sum_{n=0}^\infty \frac{1}{n!}(A+B)^n  \nonumber \\
&= \sum_{n=0}^\infty \frac{1}{n!} \sum_{k=0}^n \begin{pmatrix} n \\ k \end{pmatrix} \lbrace (A+d_B)^k 1 \rbrace B^{n-k} \\
&= \sum_{k=0}^\infty \sum_{n=k}^\infty \frac{1}{k!(n-k)!} \lbrace (A+d_B)^k 1 \rbrace B^{n-k} \nonumber \\
e^{A+B} &= [e^{A+d_B} 1 ]e^B \nonumber
\end{align}
\end{proof}

By splitting of the essential non-commutative part we get 

\subsection*{Corollary 3}
\begin{equation}
e^{A+B} = e^Ae^B + \sum_{n=0}^\infty \frac{1}{k!} D_k e^B
\end{equation}

This is different from the Campell-Baker-Hausdorff formula.

\section{Application of Theorem 1 for}
\begin{equation}
d_{B}A = hA^{2}
\end{equation}

\subsection*{Definition 3}
For $h$ a scalar and $n$ an integer we introduce
\begin{equation}
\gamma_{n}(h) = [1+h][1+2h] \cdots [1+(n-1)h], \hspace{.25cm} \gamma_{0}(h)=1
\end{equation} 

\subsection*{Lemma 1}
The following properties hold
\begin{enumerate}
\item \begin{equation*}
\gamma_{1}(h) = 1, \gamma_{n}(0) = 1, \gamma_{n}(1) = n!
\end{equation*}
\item 
\begin{equation*}
\gamma_{k+1}(h)=(1+kh)\gamma_{k}(h)
\end{equation*}
\end{enumerate}

\begin{proof}
Direct verification \\
\hfill \break
\end{proof}
Now, from Corollary 1 (14)
\begin{align*}
(A+B)^{n} &= M_{n}+ \sum_{k=2}^n \begin{pmatrix} n \\ k \end{pmatrix} D_{k}B^{n-k} \\
D_{k} &=d_{B}A^{k-1} + (A+d_{B})D_{k-1}, \hspace{.25cm} D_2=d_{B}A
\end{align*}
we find

\subsection*{Lemma 2}
\begin{enumerate}
\item \begin{equation*}
d_{B}A^{k} = khA^{k+1}
\end{equation*}
\item \begin{equation*}
D_{k} = \lbrace \gamma_{k}(h)-1 \rbrace A^{k}
\end{equation*}
\end{enumerate}

\begin{proof} \hfill \break
\begin{enumerate}
\item \begin{align*}
&d_{B}A = hA^{2} \\
&\text{Since } d_{B} \text{ is a derivation we have by induction} \\
&d_{B}A^{k} = (d_{B}A^{k-1}) A+A^{k-1}(d_{B}A) \\
&           = (k-1)h A^{k+1} + A^{k-1}hA^{2} = khA^{k+1}
\end{align*}

\item 
By induction and $D_{2} = hA^{2}$, we find \\

\begin{align*}
D_{k} &= d_{B}A^{k-1} + (A+d_{B}) \lbrace \gamma_{k-1} (h)-1 \rbrace A^{k-1} \\
&= d_{B}A^{k-1} + \lbrace \gamma_{k-1} (h)-1 \rbrace A^{k} + \gamma_{k-1}(h)d_{B}A^{k-1} - d_{B}A^{k-1} \\
&=\lbrace \gamma_{k-1} (h)-1 \rbrace A^{k} + \gamma_{k-1}(h)(k-1)hA^{k} \\
&=\lbrace [1+(k-1)h] \gamma_{k-1}(h) - 1 \rbrace A^{k} \\
D_{k} &= \lbrace \gamma_{k}(h)-1 \rbrace A^{k}
\end{align*}
\end{enumerate}
\end{proof}

Now
\begin{align*}
(A+B)^{n} &= M_{n} + \sum_{k=2}^n \begin{pmatrix} n \\ k \end{pmatrix} D_{k}B^{n-k} \\
&= \sum_{k=0}^n \begin{pmatrix} n \\ k \end{pmatrix} A^{k}B^{n-k} + \sum_{k=2}^n \begin{pmatrix} n \\ k \end{pmatrix} \lbrace \gamma_{k}(h)-1 \rbrace A^{k}B^{n-k} \\
&= B^{n} + \begin{pmatrix} n \\ 1 \end{pmatrix} AB^{n-1} + \sum_{k=2}^n \begin{pmatrix} n \\ k \end{pmatrix} \gamma_{k}(h) A^{k}B^{n-k} 
\end{align*}
Finally,
\begin{equation}
(A+B)^{n} = \sum_{k=0}^n \begin{pmatrix} n \\ k \end{pmatrix} \gamma_{k}(h) A^{k}B^{n-k}
\end{equation}

The result can also be found in \cite{4} \\
\hfill \break
Note: For $h=1$, i.e. $d_{B}A=A^{2}$, we find \\
\begin{align*}
(A+B)^{n} = &\sum_{k=0}^n \begin{pmatrix} n \\ k \end{pmatrix} k! A^{k}B^{n-k}
\end{align*}
\begin{align}
(A+B)^{n} = &\sum_{k=0}^n \frac{n!}{(n-k)!} A^{k}B^{n-k}
\end{align}

Also, if on the vector space of infinitely often differentiable function on $\mathbb{R}$ we introduce the operators

\begin{equation}
A=x, \hspace{.25cm} B=x^{2} \frac{d}{dx}
\end{equation}

we have $d_{B}A=A^{2}$. Thus the representation (21) applies.

\section{The Second Non-Commutative Binomial Theorem}

Let $A$ and $B$ be in \textswab{A}. With \\

\begin{equation}
M_{n} = \sum_{k=0}^n \begin{pmatrix} n \\ k \end{pmatrix} A^k B^{n-k} \in \textswab{A}
\end{equation}

we have

\subsection*{Lemma 3}

\begin{enumerate}

\item \begin{equation} 
M_0 = 1, M_1 = A+B 
\end{equation}

\item \begin{equation} 
M_1M_n = M_{n+1} + d_BM_n
\end{equation}

\end{enumerate}

\begin{proof} \hfill \break
\begin{enumerate}
\item Obvious
\item \begin{align*}
M_1M_n &= (A+B) \sum_{k=0}^n \begin{pmatrix} n \\ k \end{pmatrix} A^k B^{n-k} \\
&= \sum_{k=0}^n \begin{pmatrix} n \\ k \end{pmatrix} A^{k+1} B^{n-k} +  \sum_{k=0}^n \begin{pmatrix} n \\ k \end{pmatrix} BA^{k}B^{n-k} \\
&= \sum_{k=0}^n \begin{pmatrix} n \\ k \end{pmatrix} A^{k+1} B^{n-k} +  \sum_{k=0}^n \begin{pmatrix} n \\ k \end{pmatrix} \lbrace d_{B}A^{k} + A^{k}B \rbrace B^{n-k} \\
&= \sum_{s=1}^{n+1} \begin{pmatrix} n \\ s-1 \end{pmatrix} A^{s} B^{n+1-s} +  \sum_{k=0}^n \begin{pmatrix} n \\ k \end{pmatrix} A^{k}B^{n+1-k} + \sum_{k=0}^n \begin{pmatrix} n \\ k \end{pmatrix} \lbrace d_{B}A^{k} \rbrace B^{n-k} \\
&= A^{n+1} + B^{n+1} + \sum_{k=1}^{n} \left[ \begin{pmatrix} n \\ k-1 \end{pmatrix} + \begin{pmatrix} n \\ k \end{pmatrix} \right] A^{k} B^{n+1-k} + d_{B} \sum_{k=0}^n \begin{pmatrix} n \\ k \end{pmatrix} A^{k}B^{n-k} \\
&= A^{n+1} + B^{n+1} + \sum_{k=1}^{n} \begin{pmatrix} n+1 \\ k \end{pmatrix} A^{k} B^{n+1-k} + d_{B} M_{n} \\
M_1M_n &= M_{n+1} + d_{B}M_{n}
\end{align*}
\end{enumerate}

\end{proof}

\subsection*{Lemma 4}
\begin{equation}
M_{1}^{n} = M_{n} + \sum_{k=0}^{n-2} M_{1}^{k} d_{B}M_{n-1-k}
\end{equation}

\begin{proof}
This is true for $n=2$,
\begin{equation*}
M_{1}^{2} = M_{1}M_{1} = M_{2} + d_{B}M_{1}
\end{equation*}

Now by induction \\
\begin{align*}
M_{1}^{n-1} &= M_{n-1} + \sum_{k=0}^{n-3} M_{1}^{k} d_{B}M_{n-2-k} \\
M_{1}^{n} = M_{1}M_{1}^{n-1} &= M_{1}M_{n-1} + \sum_{k=0}^{n-3} M_{1}^{k+1} d_{B}M_{n-2-k} \\
&= M_{n} + d_{B}M_{n-1} + \sum_{s=1}^{n-2} M_{1}^{s} d_{B}M_{n-1-s} \\
M_{1}^{n} &= M_{n} + \sum_{k=0}^{n-2} M_{1}^{k} d_{B}M_{n-1-k} \\
\end{align*}
\end{proof}

\section*{Theorem 2}
\begin{equation}
(A + B)^{n} = M_{n} + \sum_{k=0}^{n-2} (A+B)^{k} d_{B}M_{n-1-k}
\end{equation}
\begin{proof}
This is lemma 4 with $M_1 = A+B$

\section{Application of Theorem 2 for the case}
\begin{equation}
d_{B}A = d_{B}M_1 = C, \text{ and } d_{A}C = d_{B}C = 0
\end{equation}
Then
\begin{equation}
d_{B}A^{k} = kCA^{k-1}, d_{B}M_{n} = nCM_{n-1}
\end{equation}
and
\begin{equation*}
(A+B)^n = M_{n} + \sum_{k=0}^{n-2} (n-1-k)CM_{1}^{k}M_{n-2-k}
\end{equation*}
Ansatz
\begin{equation}
(A+B)^n = \sum_{k=0}^{[\frac{n}{2}]} M_{n-2k}A_{n,k}
\end{equation}
with 
\begin{equation}
A_{n,0} = 1 
\end{equation}
and $A_{n,k}$ commuting with $A$ and $B$. \\
$[\frac{n}{2}]$ denotes the greatest integer less than $\frac{n}{2}$. \\
\hfill \break
From
\begin{equation*}
(A+B)^{n+1} = M_{1}(A+B)^{n}
\end{equation*}
we have
\begin{equation*}
\sum_{k=0}^{[\frac{n+1}{2}]} M_{n+1-2k}A_{n+1,k} = M_{1} \sum_{k=0}^{[\frac{n}{2}]} M_{n-2k}A_{n,k}
\end{equation*}
or
\begin{equation*}
M_{n+1} + \sum_{k=1}^{[\frac{n+1}{2}]} M_{n+1-2k}A_{n+1,k} = M_{1} \lbrace M_{n} + \sum_{k=1}^{[\frac{n}{2}]} M_{n-2k}A_{n,k} \rbrace
\end{equation*}
From (25) and (23) we find
\begin{equation*}
M_{1}M_{n} = M_{n+1} +nCM_{n-1}
\end{equation*}
resulting in 
\begin{equation}
\sum_{k=1}^{[\frac{n+1}{2}]} M_{n+1-2k}A_{n+1,k} = nCM_{n-1} + \sum_{k=1}^{[\frac{n}{2}]} M_{n+1-2k}A_{n,k} + \sum_{k=1}^{[\frac{n}{2}]} (n-2k)CM_{n-1-2k}A_{n,k}
\end{equation}
\hfill \break
For $n$ even, $n=2N$, (32) reads
\begin{equation*}
\sum_{k=1}^{N} M_{2N+1-2k}A_{2N+1,k} = 2NCM_{2N-1} + \sum_{k=1}^{N} M_{2N+1-2k}A_{2N,k} + \sum_{k=1}^{N-1} (2N-2k)CM_{2N-1-2k}A_{2N,k}
\end{equation*}  \\
or
\begin{align*}
M_{2N-1} A_{2N+1,1} + \sum_{k=2}^{N} M_{2N+1-2k}A_{2N+1,k} = 2NCM_{2N-1} + M_{2N-1}A_{2N,1} + \sum_{k=2}^{N} M_{2N+1-2k}A_{2N,k} \\ 
+ \sum_{k=2}^{N} M_{2N+1-2k} (2N+2-2k) A_{2N,k-1}
\end{align*}
Comparing coefficients gives the recurrence relation 
\begin{equation*}
A_{2N+1,k} = A_{2N,k} + (2N+2-2k)CA_{2N,k-1}
\end{equation*}
or
\begin{equation}
A_{n+1,k} = A_{n,k} + (n+2-2k)CA_{n,k-1}, k \geq 1
\end{equation}
Note, that for $n$ odd, $n=2N+1$, we get the same relation
\end{proof}

\subsection*{Lemma 5}
The recurrence relation (33) with $A_{n,0}=1$ has the solution
\begin{equation}
A_{n,k} = \frac{n!}{(n-2k)!k!2^k} C^{k}
\end{equation}

and (30) becomes 
\begin{equation}
(A+B)^{n} = \sum_{k=0}^{[\frac{n}{2}]} M_{n-2k} \frac{n!}{(n-2k)!k!2^{k}}C^{k}
\end{equation}

\begin{proof} by direct verification \\
This result can also be found in \cite{5} \\
\hfill \break
Note: On the vector space of infinitely often differentiable function on $\mathbb{R}$ we introduce the operators
\begin{equation}
A=x, B=\lambda\frac{d}{dx}, \text{where } \lambda \text{ is a scalar.}
\end{equation}
Then $d_{B}A=\lambda$, or $C=\lambda1$. Thus the above representation (35) applies.

\end{proof}

In particular \\
\begin{equation*}
(x+\lambda\frac{d}{dx})^{n} = \sum_{k=0}^{[\frac{n}{2}]} M_{n-2k} \frac{n!}{(n-2k)!k!2^{k}}\lambda^{k}
\end{equation*}

where \\

\begin{equation*}
M_n = \sum_{r=0}^n \begin{pmatrix} n \\ r \end{pmatrix} x^r \frac{d^{n-r}}{dx^{n-r}}, \hspace{0.5cm} M_n1 = x^n
\end{equation*}

resulting in \\

\begin{equation}
(x+\lambda \frac{d}{dx})^n 1 = \sum_{k=0}^{[\frac{n}{2}]} x^{n-2k} \frac{n!}{(n-2k)!k!2^{k}}\lambda^{k}
\end{equation}

For $\lambda= -1$, we get \\

\begin{equation}
(x-\frac{d}{dx})^n 1 = n! \sum_{k=0}^{[\frac{n}{2}]} (-1)^{k} \frac{x^{n-2k}}{(n-2k)!k!2^{k}}
\end{equation}

The right-hand side are the Hermite polynomials. \\

Thus \\

\begin{equation}
He_n(x) = (x-\frac{d}{dx})^n 1
\end{equation}

\noindent\textit{Department of Physics, University of Colorado Boulder, Boulder, CO 80309\\
Walter.Wyss@Colorado.EDU}

\end{document}